\documentclass{article}

\usepackage{amsmath}
\usepackage{amsfonts}
\usepackage{amssymb}
\usepackage{amsthm}
\usepackage[all]{xy}

\theoremstyle{plain}

\newtheorem{theorem}{Theorem}[section]

\newtheorem{proposition}[theorem]{Proposition}

\newtheorem{example}[theorem]{Example}
\newtheorem{examples}[theorem]{Examples}
\theoremstyle{definition}
\newtheorem{definition}[theorem]{Definition}

\newtheorem{remark}[theorem]{Remark}

\def\To{\mathrel{\xy \ar@{=>} (200,0) \endxy}}

\begin{document}

\title{On closure operators and reflections \\ in Goursat categories}

\author{Francis Borceux\footnote{D\'epartement de Math\'ematique, Universit\'e de Louvain,
2 chemin du Cyclotron,
1348 Louvain-la-Neuve, Belgique, e-mail: \texttt{francis.borceux@uclouvain.be}},\, Marino Gran\footnote{Marino Gran, Lab.\ Math.\ Pures et Appliqu\'ees, Universit\'e du Littoral C\^ote d'Opale, 50,~Rue F.~Buisson, 62228~Calais, France. Email:~\texttt{gran@lmpa.univ-littoral.fr}}\, and Sandra Mantovani\footnote{Dipartimento di Matematica, Universit\`a degli Studi di Milano,
 Via C. Saldini 50,
 20133 Milano, Italy, e-mail: \texttt{sandra.mantovani@mat.unimi.it}
 }
}
\date{\emph{Dedicated to the memory of Fabio Rossi}}
\maketitle

\begin{abstract}
\noindent By defining a closure operator on effective equivalence relations in a regular category $\mathcal C$,
it is possible to establish a bijective correspondence between these closure operators and the regular epireflective
subcategories $\mathcal L$ of $\mathcal C$. When $\mathcal C$ is an exact Goursat category this correspondence
restricts to a bijection between the Birkhoff closure operators on effective equivalence relations and the Birkhoff
subcategories of $\mathcal C$. In this case it is possible to provide an explicit description of the closure,
and to characterise the congruence distributive Goursat categories.

\smallskip

\end{abstract}

\section{Introduction}
A classical result in the theory of abelian categories establishes a bijection between the hereditary torsion
theories and the universal closure operators. This correspondence has been recently extended to the
non-additive context of homological categories in \cite{BG}; more precisely,
the correspondence has been shown to be nothing but the restriction of a
general bijection between the (regular-)\\epireflective subcategories and the so-called homological closure operators.

It is the aim of the present note to discuss a further extension of this last
result in any regular category, and to give an explicit description of the closure operator
in various contexts which are relevant in universal algebra and in topological algebra.
When the base category $\mathcal C$ is homological it is possible to define
the closure operator on \emph{normal subobjects} (=kernels), because these
subobjects correspond to quotients, as in the familiar case of the category of groups.
Of course, in a regular category this is no longer possible: the natural thing to do in this
non-pointed situation is then to define the closure on \emph{effective equivalence relations}.
We determine in Theorem \ref{Effective} the conditions on such a closure operator that make
it correspond to an epireflective subcategory $\mathcal L$ of $\mathcal C$. As expected,
the closure arising in this way from the reflection of the regular category ${\mathbb T} (\mathbf{Top})$
of Mal'tsev topological algebras into its subcategory ${\mathbb T} (\mathbf{Haus})$ of Hausdoff topological
algebras coincides with the usual topological closure (Example \ref{topological}).

We then restrict ourselves to the so-called \emph{exact Goursat categories} \cite{CKP}: these categories have
the property that the composition of two equivalence relations $R$ and $S$ on any fixed object is $3$-permutable:
$$R \circ S \circ R = S \circ R \circ S.$$
These categories are very common in universal algebra, and they
include all Mal'tsev varieties \cite{S}: in particular groups,
quasi-groups, rings, associative algebras, Heyting algebras and
implication algebras are all examples of Goursat varieties. If the
epireflective subcategory $\mathcal L$ of an exact Goursat
category $\mathcal C $ is stable in $\mathcal C $ under regular
quotients, i.e. if $\mathcal L$ is a Birkhoff subcategory of
$\mathcal C$, then the closure is preserved by regular images, so
that $f(\overline{S})=\overline{f(S)}$ for any equivalence
relation $S$ on $X$ and any regular epimorphism $f \colon X
\rightarrow Y$. Furthermore, the closure of $S$ is given by the
formula
$$\overline{S}= \overline{\Delta_X} \circ S \circ \overline{\Delta_X},$$
where $\Delta_X$ is the discrete equivalence relation on $X$ (Proposition \ref{formula}). This means that the knowledge
of the closure of $\Delta_X$ suffices already to determine the closure of any equivalence relation $S$ on $X$. We can use
this formula to get the closure determined by the reflection of the exact category ${\mathbb T} (\mathbf{HComp})$ of
Mal'tsev compact Hausdorff algebras into its subcategory ${\mathbb T} (\mathbf{Profin})$ of profinite topological
algebras (Proposition \ref{profinite}). The closure of an equivalence relation $S$ on $A$ is given here by
$\overline{S}= S \circ R_A,$ where $R_A$ is the congruence on $A$ that identifies two points when they are
in the same connected component.

Finally, by using a recent result of Bourn in \cite{Bourn}, it is possible to obtain a characterisation of the exact Goursat categories that are congruence distributive, in terms of a property of the
closure operator associated with any Birkhoff subcategory, namely
the preservation under regular images of
 the closure of an intersection : $$f(\overline{R \wedge S})=\overline{f(R)} \wedge \overline{f(S)}.$$

\section{Reflections in regular categories}

In this section $\mathcal C$ will denote a finitely complete regular category.
By an {\em epireflective subcategory} $\mathcal L$ of $\mathcal C$ we shall mean a full replete reflective subcategory
$$
\xymatrix@=30pt{
{\mathcal L} \ar@<-1ex>[r]_{\iota} & {\mathcal C}
\ar@<-1ex>[l]^{\perp}_{\lambda}}
$$
with the property that every component $\eta_A \colon A \rightarrow \iota \lambda(A)$ of the unit of the adjunction is a regular epimorphism; this implies in particular that $\mathcal L$ will always be closed in $\mathcal C$ under subobjects.

Recall that an equivalence relation $(S, p_1,p_2)$ on an object $X$ is said to be effective, when it is the kernel pair $R[f]$ of a morphism $f \colon X \rightarrow Y$. In a regular category this is the case if and only if $(S, p_1,p_2)$ is the kernel pair of the coequaliser of the two projections $p_1 \colon S \rightarrow X$ and $p_2 \colon S \rightarrow X$. By an exact category \cite{Ba} we shall always mean a finitely complete regular category with the property that equivalence relations are effective.
We would like to characterise the epireflective subcategories of a given regular category $\mathcal C$ in terms of a special kind of closure operator on effective equivalence relations. For this, the following definitions will be needed:

\begin{definition}
An {\em idempotent closure operator $\overline{(\, \,)}$ on effective equivalence relations} consists in giving for every effective equivalence relation $\xymatrix{S \ar[r]^s&  X \times X}$ another effective equivalence relation $\xymatrix{\overline{S} \ar[r]^{\overline{s}} &  X \times X}$ called the closure of $S$ in $X \times X$. This assignment has to satisfy the following properties, where $S$ and $T$ are effective equivalence relations on $X$, $\xymatrix{ Y \ar[r]^{f}&  X}$ is an arrow in $\mathcal L$ and $f^{-1} (S)$ is the inverse image of $S$ along $f$:
\begin{enumerate}
\item $S \subseteq \overline{S}$
\item $S \subseteq T$ implies $\overline{S} \subseteq \overline{T}$
\item $\overline {f^{-1} (S)} \subseteq f^{-1} (\overline{S})$
\item $\overline{\overline{S}} = \overline{S}$
\end{enumerate}
\end{definition}
We shall write $i_S \colon S \rightarrow \overline{S}$ for the canonical inclusion of $S$ in $\overline{S}$.
An effective equivalence relation $\xymatrix{S \ar[r]^s&  X \times X}$ is {\em closed} when $S = \overline{S}$, and {\em dense} when $\overline{S}=X \times X$.
\begin{definition}
An idempotent closure operator on effective equivalence relations will be called an {\em effective closure operator} (on effective equivalence relations) if it also satisfies the following axiom: \\
$(5)$ for any regular epimorphism $\xymatrix{ Y \ar[r]^{f}&  X}$ one has that $\overline {f^{-1} (S)} = f^{-1} (\overline{S}).$
\end{definition}

The following result can be considered as a ``non-pointed'' version of Theorem $2.4$ in \cite{BG}:
\begin{theorem}\label{Effective}
Let $\mathcal C$ be a regular category.
There is a bijection between the epireflective subcategories of $\mathcal C$ and the effective closure operators.
\end{theorem}
\begin{proof}
Let us begin with an epireflective subcategory $\mathcal L$ of $\mathcal C$ with reflector $\lambda \colon {\mathcal C} \rightarrow {\mathcal L}$. For any object $X$ in $\mathcal C$ there is a canonical exact fork
$$ \xymatrix@=20pt{R[\eta_X]
\ar@<-1ex>[r]_{p_2}
\ar@<1ex>[r]^{p_1}
 & X  \ar[r]^{\eta_X}   &  \iota \lambda(X)}
$$
with $(R[\eta_X],p_1,p_2)$ the kernel pair of the $X$-component $\eta_X$ of the unit of the adjunction.

The closure $\overline{S}$ of an effective equivalence relation $S$ is defined as the inverse image of the equivalence relation $R[\eta_{X/S}]$ along the quotient $q$ of $X$ by $S$, i.e. $\overline S = q^{-1} (R[\eta_{X/S}])$:
$$ \xymatrix@=40pt{
 & \overline{S} \ar@<1ex>[d] \ar@<-1ex>[d]^{} \ar[r]^{}   & R[\eta_{X/S}] \ar@<1ex>[d]^{p_2} \ar@<-1ex>[d]_{p_1}  \\
 S \ar@<1ex>[r] \ar@<-1ex>[r] \ar@{.>}[ru]^{i_S}
 & X  \ar[r]_{q}   & X/S
 }
$$
It is clear that the equivalence relation $\overline{S}$ could be equivalently defined as the kernel pair of $\eta_{X/S} \circ q : X \rightarrow \iota \lambda (X/S)$, so that $\overline{S}= R[\eta_{X/S} \circ q]$.
It follows that there is an inclusion $i_S \colon S \rightarrow \overline{S}$, showing the validity of axiom (1).
Remark that, in particular, the kernel pair $R[\eta_X]$ of the reflection $\eta_X \colon X \rightarrow \iota \lambda (X) $ is exactly the closure
$ {\overline{\Delta}_X}$ of the discrete equivalence relation ${\Delta}_X$ of $X$. Since $\eta_{X/S} \circ q$ is the coequaliser of its kernel pair $\overline{S}$, the axiom $(4)$ is also satisfied. The axiom $(2)$ is a consequence of the fact that $\lambda$ is a functor.

Consider then the diagram
$$
 \xymatrix@=30pt{f^{-1}(S) \ar[d]_g
\ar@<-1ex>[r] \ar@<1ex>[r]
 & Y \ar[d]^{f} \ar[r]^(.4){q'}   & Y/f^{-1}(S) \ar[d]^{h} \ar@{.>}[r]& \lambda(Y/f^{-1}(S))\ar[d]^{\lambda(h)} \\
 S \ar@<-1ex>[r] \ar@<1ex>[r]
 & X  \ar[r]_q   & X/S \ar@{.>}[r]& \lambda(X/S)
 }
$$
where $q$ and $q'$ are the canonical quotients. The functoriality of $\lambda$ and the construction of the closure operator on effective equivalence relations give a unique arrow $\tau$ such that the square
$$ \xymatrix@=25pt{\overline{f^{-1}(S)} \ar@{}[dr]|{(\star)}\ar[r]^(.6){\overline{f^{-1}(s)}}
\ar@{.>}[d]_{\tau}
 & {Y \times Y}  \ar[d]^{f}   \\
 \overline{S} \ar[r]_{\overline{s}}
 & {X \times X}
 }
$$
commutes. This implies that $\overline {f^{-1} (S)} \subseteq f^{-1} (\overline{S})$, proving axiom $(3)$.

Remark that, by construction of the inverse image of an equivalence relation, one knows that the induced arrow $h \colon Y/f^{-1}(S) \rightarrow X/S$ in the diagram above is always a monomorphism.
When, moreover, the arrow $f$ is assumed to be a regular epimorphism, the arrow $h \colon Y/f^{-1}(S) \rightarrow X/S$ actually is an isomorphism. This easily implies that, in this case, the square $(\star)$ here above is a pullback, so that axiom $(5)$ holds true.

Conversely, given an effective closure operator on effective equivalence relations $\overline{(\, \,)}$, we are going to prove that the full replete subcategory $\mathcal L$ of $\mathcal C$ defined by
$$X \in {\mathcal L} \quad {\rm if \, and \,  only \, if} \quad \Delta_X {\rm \, is \, closed} $$
is epireflective in $\mathcal C$. Given an object $X$ in $\mathcal C$, in order
to define the left adjoint $\lambda \colon \mathcal C \rightarrow \mathcal L$, we consider the closure $\overline{\Delta}_X$ of the discrete equivalence relation $\Delta_X$ on $X$, and the canonical quotient ${X}/\overline{\Delta}_X$
$$ \xymatrix@=20pt{\overline{\Delta}_X \ar@<-1ex>[r]_{p_2}
\ar@<1ex>[r]^{p_1}
 & X  \ar[r]^(.4){\eta_X} & {X}/\overline{\Delta}_X ,}
$$
where $\eta_X$ is the coequaliser of $p_1$ and $p_2$. We define the reflector on objects by setting $\lambda(X)= {X}/\overline{\Delta}_X$, for any $X$ in $\mathcal C$. In order to see that the object ${X}/\overline{\Delta}_X$ belongs to $\mathcal L$, we first observe that, by axiom $(5)$, $\eta_X^{-1} (i_{\Delta_{{X}/\overline{\Delta}_X }})= i_{\overline{\Delta}_X}$. But the arrow $i_{\overline{\Delta}_X} \colon \overline{\Delta}_X \rightarrow \overline{\overline{\Delta}_X}$ is an isomorphism by axiom $(4)$, from which it easily follows that the canonical inclusion $i_{\Delta_{{X}/\overline{\Delta}_X }} \colon \Delta_{{X}/\overline{\Delta}_X } \rightarrow \overline{\Delta_{{X}/\overline{\Delta}_X }}$ is a regular epimorphism, thus an isomorphism, as desired.

Let us then show that $\eta_X \colon X \rightarrow X/\overline{\Delta}_X$ has the universal property with respect to the full subcategory $\mathcal L$. For this, let $f \colon X \rightarrow Y$ be any arrow with $Y$ in $\mathcal L$, so that $\overline{\Delta}_Y = \Delta_Y$:
$$ \xymatrix@=20pt{
  \overline{\Delta}_X  \ar@<-1ex>[r]_{p_2}
\ar@<1ex>[r]^{p_1}  & X \ar[d]^{f} \ar[r]^{\eta_X}& X/{\overline{\Delta}_X}\\
 \Delta_Y = \overline{\Delta}_Y     \ar@<-1ex>[r]_(.6){p_2}
\ar@<1ex>[r]^(.6){p_1}    & Y &
 }
$$
Now, since $\Delta_X \subseteq R[f]= f^{-1}(\Delta_Y)$, by axioms $(2)$ and $(3)$ it follows that
$$\overline{\Delta}_X \subseteq \overline{f^{-1}(\Delta_Y)} \subseteq f^{-1}(\overline{\Delta}_Y)= f^{-1}({\Delta}_Y)=R[f].$$
Accordingly, the universal property of the coequaliser $\eta_X$ gives a unique arrow $g \colon X/{\overline{\Delta}_X} \rightarrow Y$ such that $g\circ \eta_X = f$.
 The functor $\lambda \colon \mathcal C \rightarrow \mathcal L$ is then the left adjoint of the forgetful functor $\iota \colon \mathcal L \rightarrow \mathcal C$, and $\mathcal L$ is epireflective in $\mathcal C$.

Finally, let us check that these constructions determine a bijection between epireflective subcategories of $\mathcal C$ and effective closure operators. If $\overline{(\, \,)}$ is a special closure operator, $\lambda$ is the functor defined as above, and $\widetilde{(~)}$ is the closure operator on effective equivalence relations associated with $\lambda$, we have to prove that $\overline{(\, \,)}= \widetilde{(~)}$. For this, it suffices to consider the following diagram:
$$ \xymatrix@=20pt{
 & {\widetilde{S}} \ar@<1ex>[d] \ar@<-1ex>[d] \ar[r]^{}   & \overline{\Delta}_{X/S} \ar@<-1ex>[d]^{} \ar@<1ex>[d]^{} \\
 S \ar@<1ex>[r]_{} \ar@<-1ex>[r]_{} \ar[ru]^{}
 & X  \ar[r]_{q}   & X/S }
$$
By axiom $(5)$ of the original closure operator $\overline{(\, \,)}$ it follows that $$\widetilde{S}=q^{-1}(\overline{\Delta}_{X/S} ) = \overline{q^{-1}(\Delta_{X/S})} = \overline{S}.$$
Conversely, given an epireflection $\lambda \colon \mathcal C \rightarrow \mathcal L$, the associated closure operator $\overline{(\, \,)}$ and the epireflection $\overline{\lambda} \colon \mathcal C \rightarrow \overline{\mathcal L}$ associated with $\overline{(\, \,)}$, we have to prove that $\mathcal L = \overline{\mathcal L}$. This follows easily from the fact that the quotient $X / \overline{\Delta}_X$ is exactly the reflection of $X$ in $\overline{\mathcal L}$.
\end{proof}
\begin{example}\label{topological}
\rm{Consider any Mal'tsev theory $\mathbb T$, i.e. is any algebraic theory containing a ternary term $p(x,y,z)$ satisfying the identities $p(x,x,y)=y$ and $p(x,y,y)=x$ \cite{S}. Any algebraic theory containing a group operation is in particular a Mal'tsev theory, since it suffices to set $p(x,y,z)= x \cdot y^{-1} \cdot z$ to get a Mal'tsev operation; however, there are also other interesting examples, such as the theories of quasi-groups and of Heyting algebras \cite{Jo}. The category ${\mathbb T} (\mathbf{Top)}$ of topological models of such a theory is a regular Mal'tsev category in the sense of Carboni, Lambek and Pedicchio \cite{CLP}. Regular epimorphisms in this category are given by the open surjective homomorphisms  \cite{Hus, JP}, as in the classical case of the category ${Grp} (\mathbf{Top)}$ of topological groups. An object in ${\mathbb T} (\mathbf{Top)}$ is called a \emph{topological Mal'tsev algebra}.

The category ${\mathbb T} (\mathbf{Top})$ of topological Mal'tsev algebras is reflective in the category $\mathbb T (\mathbf{Haus})$
of Hausdorff Mal'tsev algebras: the reflection of an algebra $A$ is obtained by taking the quotient $\pi^{\mbox{\itshape\footnotesize Top}}
 \colon A \rightarrow \frac{A} {\overline{\Delta_A}^{\mbox{\itshape\footnotesize Top}}
}$ of $A$ by the \emph{topological closure} $\overline{\Delta_A}^{\mbox{\itshape\footnotesize Top}}
$ of the discrete equivalence relation ${\Delta_A}$ on $A$. As observed in \cite{BC}, the topological closure $\overline{\Delta_A}^{\mbox{\itshape\footnotesize Top}}
$ of $\Delta_A$ is automatically a congruence, since it is a reflexive relation in $\mathbb T (\mathbf{Top)}$, and this latter is a Mal'tsev category. We then get the epireflective subcategory
$$
\xymatrix@=30pt{
{\mathbb T} (\mathbf{Haus})
\ar@<-1ex>[r]_{\iota} & {\mathbb T (\mathbf{Top})}
\ar@<-1ex>[l]^{\perp}_{\lambda}}
$$
and we are going to show that the effective closure $\overline{(\, \,)}$ associated with this reflection as in Theorem \ref{Effective} coincides with the usual topological closure $\overline{(\, \,)}^{\mbox{\itshape\footnotesize Top}}$. By construction of the effective closure operator one clearly has that $\overline{\Delta_A} =\overline{\Delta_A}^{\mbox{\itshape\footnotesize Top}}
$, and this implies that, for any effective equivalence relation $S$ on $A$,
$\overline{S}^{\mbox{\itshape\footnotesize Top}}
 \subseteq \overline{S}$.
\\
On the other hand, the fact that $\overline{S}^{\mbox{\itshape\footnotesize Top}}
$ is saturated with respect to the open surjective homomorphism $$q^{\mbox{\itshape\footnotesize Top}}
 \times q^{\mbox{\itshape\footnotesize Top}}
 \colon A \times A \rightarrow \frac{A} {\overline{S}^{\mbox{\itshape\footnotesize Top}}
} \times \frac{A} {\overline{S}^{\mbox{\itshape\footnotesize Top}}
}$$
 implies that $\frac{A} {\overline{S}^{\mbox{\itshape\footnotesize Top}}
}$
is a Hausdorff algebra, as its diagonal $\Delta$ is then closed. This yields an arrow $\alpha \colon i \lambda(\frac{A}{S}) \rightarrow \frac{A} {\overline{S}^{\mbox{\itshape\footnotesize Top}}
}$ with $\alpha \cdot \eta_{\frac{A}{S}} \cdot q = q^{\mbox{\itshape\footnotesize Top}}
$, giving the inclusion
$\overline{S} \subseteq \overline{S}^{\mbox{\itshape\footnotesize Top}}
$. It follows that $\overline{S} = \overline{S}^{\mbox{\itshape\footnotesize Top}}
$.}
\end{example}
\section{Birkhoff subcategories of a Goursat category}
In this section we restrict ourselves to the context of exact Goursat categories. This is the natural context where Birkhoff subcategories can be characterised by a special kind of closure operators on equivalence relations.\\
We begin by recalling a few well known definitions:
\begin{definition} \cite{JK}
Let $\mathcal C$ be an exact category, and let $\mathcal L$ be an epireflective subcategory of $\mathcal C$. One calls $\mathcal L$ a \emph{Birkhoff subcategory} of $\mathcal C$ when it is closed under regular quotients in $\mathcal C$: if $f \colon X \rightarrow Y$ is a regular epimorphism with $X$ in $\mathcal L$, then $Y$ belongs to $\mathcal L$ as well.
\end{definition}
Of course, when $\mathcal C$ is a variety of universal algebras, $\mathcal L$ is a Birkhoff subcategory of $\mathcal C$ exactly when it is a subvariety.
\begin{definition}
An exact category $\mathcal C$ is a \emph{Goursat category} \cite{CKP} when any pair of equivalence relations $R$, $S$ on the same object $X$ in $\mathcal C$ satisfies the condition $$R \circ S \circ R = S \circ R \circ S,$$ where the symbol $\circ$ denotes here the usual relational composition in a regular category \cite{Ba}.
\end{definition}
\begin{example}
\rm{The notion of Goursat category arises from universal algebra. Indeed, Goursat varieties are well known under the name of $3$-permutable varieties. These varieties are characterised by a Mal'tsev condition as follows: a variety is $3$-permutable if and only if there exist two ternary terms $p$ and $q$ satifying $p(x,y,y)=x$, $q(x,x,y)=y$ and $p(x,x,y)=q(x,y,y)$ \cite{HM}.
Of course, any Mal'tsev category
\cite{CLP} is in particular a Goursat category, since the $2$-permutability of the composition of equivalence relations is a stronger condition than the one of $3$-permutability. A genuine example of a Goursat category which is not a Mal'tsev category is provided by the variety of implication algebras (see \cite{Gumm}). These algebras are equipped with a single binary operation $\triangleright$ satisfying the identities $(x \triangleright y) \triangleright y = (y \triangleright x) \triangleright x$, $(x \triangleright y) \triangleright x =x $ and $x \triangleright (y \triangleright z)= y \triangleright ( x \triangleright z)$.}
\end{example}
\begin{definition}
A \emph{Birkhoff closure operator} on effective equivalence relations $\overline{(\, \,)}$ is an effective closure operator satisfying the following additional property $(6)$: for any regular epi
$f \colon X \rightarrow Y$ and any equivalence relation $S \rightarrow X \times X$ one has that \\
$$(6 )\qquad \overline{f(S)}= f(\overline{S}).$$
\end{definition}
\begin{remark}
It was shown in \cite{CKP} that a regular category has the Goursat property if and only if the regular image $f(S)$ of an equivalence relation $S$ is an equivalence relation.
Therefore the axiom $(6)$ here above only makes sense in an exact Goursat category: in order to consider the closure $\overline{f(S)}$ of $f(S)$ we need to know that $f(S)$ is an effective equivalence relation, and this is always the case exactly under the Goursat assumption.
\end{remark}
\begin{proposition}
Let $\mathcal C$ be an exact Goursat category.
There is a bijection between the Birkhoff subcategories of $\mathcal C$ and the Birkhoff closure operators.
\end{proposition}
\begin{proof}
We must show that an effective closure operator on equivalence relations $\overline{(\, \,)}$ satisfies axiom $(6)$ if and only if the corresponding epireflective subcategory $\mathcal L$ is closed in $\mathcal C$ under regular quotients.

On the one hand, let axiom $(6)$ hold, and let $f\colon X \rightarrow Y$ be a regular epi with $X\in \mathcal L$. One then gets that
$$\overline{\Delta}_Y = \overline{f(\Delta_X)}= f(\overline{\Delta}_X) = f(\Delta_X) = \Delta_Y,$$
so that $\Delta_Y$ is closed, and $Y \in \mathcal L$.

On the other hand, let us assume that the epireflective subcategory $\mathcal L$ is closed in $\mathcal C$ under regular quotients. Given a regular epi $f \colon X \rightarrow Y$ and an equivalence relation $S$ on $X$, its regular image $f(S)$ along $f$ gives rise to the following diagram:

$$
 \xymatrix@=30pt{S \ar[d]_g
\ar@<-1ex>[r] \ar@<1ex>[r]
 & X \ar[d]^{f} \ar@{}[dr]|{(1)} \ar[r]^(.4){q'}   & X/S \ar@{}[dr]|{(2)} \ar[d]^{h} \ar@{.>}[r]^{\eta_{X/S}} & \lambda(X/S)\ar[d]^{\lambda(h)} \\
 f(S) \ar@<-1ex>[r] \ar@<1ex>[r]
 & Y  \ar[r]_(.4){q}   & Y/f(S) \ar@{.>}[r]_{\eta_{Y/f(S)}} & \lambda(Y/f(S))
 }
$$
The fact that $g$ is an epimorphism implies that the square $(1)$ is a pushout; furthermore, the square $(2)$ is also a pushout, precisely because the arrow $h$ is a regular epimorphism and the subcategory $\mathcal L$ is closed in $\mathcal C$ under regular quotients (see Proposition $3.1$ in \cite{JK}). Now, in an exact Goursat category, the fact that the rectangle $(1)+(2)$ is a pushout of regular epis implies that the induced arrow $\tilde{f} \colon R[\eta_{X/S} \circ q'] \rightarrow R[\eta_{Y/f(S)} \circ q]$ from the kernel pair $R[\eta_{X/S} \circ q'] = \overline{S}$ to the kernel pair $R[\eta_{Y/f(S)} \circ q] = \overline{f(S)}$ is a regular epi. Indeed, the fact that $(1)+(2)$ is a pushout clearly implies that $\eta_{Y/f(S)} \circ q$ is the coequaliser of the projections $p_1$ and $p_2$ of the relation $f(\overline{S})$ on $Y$; but this latter is an equivalence relation (by the Goursat property), which is effective (by exactness), so that $f(\overline{S})= R[\eta_{Y/f(S)} \circ q] = \overline{f(S)}$, as desired.
\end{proof}

\begin{proposition}\label{formula}
Let $\mathcal C$ be an exact Goursat category, $\mathcal L$ a Birkhoff subcategory of $\mathcal C$, $\overline{(\, \,)}$ the corresponding Birkhoff closure operator on equivalence relations. \\
Then, for any equivalence relation $S$ on any $X$ in $\mathcal C$, one has the following formula for the closure $\overline{S}$ of $S$:
$$\overline{S}= \overline{\Delta_X} \circ S \circ \overline{\Delta_X}=S \circ \overline{\Delta_X} \circ S. $$
\end{proposition}
\begin{proof}
By definition of the closure operator on the equivalence relations, the closure $\overline{S}$ is obtained as the kernel pair of the diagonal $\eta_{X/S} \circ q$ of the following commutative square
$$ \xymatrix@=30pt{
 X \ar[r]^{q} \ar[d]_{\eta_X} & X/S \ar[d]^{\eta_{X/S}}  \\
  \iota \lambda(X) \ar[r]_{\iota \lambda (q)}   & \iota \lambda (X/S)
 }
$$
which is a pushout by the Birkhoff assumption. This precisely means that $$\overline{S}= \overline{\Delta_X}\vee S.$$ In an exact Goursat category, the join of the two equivalence relations $\overline{\Delta_X}$ and $S$ in the modular lattice of equivalence relations on $X$ always exists, and is given exactly by the formula $$ \overline{\Delta_X}\vee S = \overline{\Delta_X} \circ S \circ \overline{\Delta_X}=S \circ \overline{\Delta_X} \circ S.$$
\end{proof}
\begin{remark}
It is clear that for an effective closure operator the property $(6)$ is also equivalent to the following (apparently) weaker property: for any regular epimorphism $f \colon X \rightarrow Y$ one has that
$$(6') \qquad f(\overline{\Delta_X})= \overline{\Delta_Y}. $$
Indeed, to check that $(6')$ implies $(6)$ it suffices to observe that
$$\overline{f(S)}= f(S) \vee \overline{\Delta_Y} = f(S) \vee f(\overline{\Delta_X}) = f(S \vee \overline{\Delta_X})= f(\overline{S} ).$$
\end{remark}
\begin{example}
\rm{
Consider again a Mal'tsev theory $\mathbb T$ and, this time, the category of compact Hausdorff Mal'tsev algebras ${\mathbb T} (\mathbf{HComp})$. ${\mathbb T} (\mathbf{HComp})$ is an exact Mal'tsev category, with regular epimorphisms given by open (and closed) surjective homomorphisms. We denote by ${\mathbb T} (\mathbf{Profin})$ the category of profinite (=compact Hausdorff totally disconnected) Mal'tsev algebras. We are now going to explain why ${\mathbb T} (\mathbf{Profin})$ is a Birkhoff subcategory of ${\mathbb T} (\mathbf{HComp})$, and describe the corresponding closure. In order to make the paper more self-contained, we repeat here the arguments of \cite{BC1, BC, BG} in the protomodular and semi-abelian cases, which are still valid in this more general context.}
\begin{proposition}\label{profinite}
Let  ${\mathbb T}$ be a Mal'tsev theory.
\begin{enumerate}
\item ${\mathbb T} (\mathbf{Profin})$ is a Birkhoff subcategory of ${\mathbb T}(\mathbf{HComp})$;
\item the effective closure of an equivalence relation $S$ on $A$ is given by $$\overline{S}= S \circ R_A,$$
    where $R_A$ is the congruence on $A$ that identifies two points when they are in the same connected component.
\end{enumerate}
\end{proposition}
\begin{proof}
$1.$
Given a compact Hausdorff Mal'tsev algebra $A$, we write $\Gamma(a)$ for the connected component of the element $a$ in $A$. Write $R_A$ for the subset $$R_A = \{(a,b) \in A \times A \mid \Gamma(a)= \Gamma(b) \} .$$
By using the fact that the Mal'tsev operation $p_A \colon A \times A \times A \rightarrow A$ is continuous, and that the continuous image of a connected space is connected, one can check that $\Gamma(p_A(a_1,a_2,a_3))= \Gamma(p_A(b_1,b_2,b_3))$ whenever $\Gamma(a_1)= \Gamma(b_1)$, $\Gamma(a_2)= \Gamma(b_2)$ and $\Gamma(a_3)= \Gamma(b_3)$, so that $R_A$ is a subalgebra of $A \times A$. Since $R_A$ is reflexive and ${\mathbb T} (\mathbf{HComp})$ is a Mal'tsev category, $R_A$ is a congruence on $A$. The reflection $\lambda \colon {\mathbb T} (\mathbf{HComp}) \rightarrow {\mathbb T} (\mathbf{Profin})$ is given, for any $A$ in ${\mathbb T} (\mathbf{HComp})$, by the quotient
$\eta_A \colon A \rightarrow \frac{A}{R_A}$. Indeed, the algebra $\frac{A}{R_A}$ is compact as continuous image of a compact space; to see that it is totally disconnected, consider, for any $a$ in $A$, the pullback
$$ \xymatrix@=30pt{
 P \ar[r]^{p_2} \ar[d]_{p_1} & \Gamma([a]) \ar[d]^{i}  \\
 A \ar[r]_{\eta_A}   & \frac{A}{R_A},
 }
$$
where $i$ is the canonical inclusion of the connected component $\Gamma([a])$ of the equivalence class $[a]$ in the quotient $\frac{A}{R_A}$. The arrow $p_2$ is an open surjection, with connected codomain and connected fibres: the $q$-reversibility of connected spaces \cite{AW} then implies that $P$ is a connected space. It follows that $P = \Gamma(a)$, and $\Gamma([a])= [a]$: the algebra $\frac{A}{R_A}$ is then totally disconnected, actually the profinite reflection of the algebra $A$.

It remains to show that the epireflective subcategory ${\mathbb T}
(\mathbf{Profin})$ is stable in ${\mathbb T} (\mathbf{HComp})$
under regular quotients. Consider an open surjective homomorphisms
$f \colon A \rightarrow B,$ with $A$ a profinite algebra. Since
$B$ is compact Hausdorff, to show that $B$ is totally disconnected
reduces to prove that for any $b$ and $b'$ in $B$, with
$b\not=b'$, there exist in $B$ two disjoint open and closed
subsets separating these two points. By the assumption of
profiniteness of $A$, there exist two disjoint open and closed
subsets $U$ and $U'$ of $A$ with the property that $f^{-1} (b)
\subset U$ and $f^{-1} (b') \subset U'$. Since $f$ is open and
closed, $f(U)$ is open and closed, as are $f^{-1} (f (U))$ and its
complement $A \setminus  f^{-1} (f (U))$. Since $U \cap f^{-1}
(b')= \emptyset$, $f^{-1} (b') \subset A\setminus f^{-1} (f (U))$
and
it follows that $f(U) $ and $f(A\setminus f^{-1} (f (U))$ are the two closed and open subsets of $B$ we were looking for. \\

$2.$ It follows from $1.$ and Proposition \ref{formula}, by taking into account the fact that
$S \vee R_A = S \circ R_A$ in the exact Mal'tsev category ${\mathbb T} (\mathbf{HComp})$.
\end{proof}

\end{example}
As a consequence of Proposition \ref{formula}, we can prove that the closure operator on equivalence relations preserves binary joins, i.e. it is an additive closure operator \cite{DT}:
\newpage
\begin{proposition}
Let $\mathcal C$ be an exact Goursat category, $\mathcal L$ a Birkhoff subcategory of $\mathcal C$, $\overline{(\, \,)}$ the corresponding Birkhoff closure operator. Then:
\begin{enumerate}
\item given a regular epimorphism $f \colon X \rightarrow Y$ and two equivalence relations $R$ and $S$ on $X$, one has that $$\overline{f(R \vee S)}= \overline{f(R)} \vee \overline{f(S)}; $$
\item for any equivalence relations $R$ and $S$ on $X$, one has $$\overline{R \vee S}= \overline{R} \vee \overline{S}.$$
\end{enumerate}
\end{proposition}
\begin{proof}
$1.$ The Goursat assumption implies that the regular images $f(R)$, $f(S)$ and $f(R\vee S)$ are (effective) equivalence relations, and that the join of two (effective) equivalence relations always exists. By property $(2)$ of the corresponding closure operator $\overline{(\, \,)}$, the equivalence relation $\overline{f(R \vee S)}$ contains both $\overline{f(R)}$ and $\overline{f(S)}$, thus  $\overline{f(R \vee S)} \geq \overline{f(R)} \vee \overline{f(S)}. $ \\
In order to check the other inclusion, one first observes that
$$f(R \vee S)=f(R) \vee f(S):$$
this easily follows from
the fact that
$f(R)$, $f(S)$ and $f(R\vee S)$ are equivalence relations, by using the construction of the join of two equivalence relations as the kernel pair of the diagonal of the pushout of the corresponding quotients.
Next, by using the description of the closure given in Proposition \ref{formula} one gets that
\begin{eqnarray*} \overline{f(R \vee S)} & = & \overline{f(R) \vee f(S)} \\
& = & \overline{\Delta}_Y \circ f(R) \circ f(S) \circ f(R) \circ \overline{\Delta}_Y \\
 & \leq &  ( \overline{\Delta}_Y \circ f(R) \circ \overline{\Delta}_Y ) \circ \overline{f(S)} \circ ( \overline{\Delta}_Y \circ f(R) \circ \overline{\Delta}_Y ) \\
 & = & \overline{f(R)} \circ \overline{f(S)}  \circ \overline{f(R)} \\
 & = &  \overline{f(R)} \vee  \overline{f(S).}
  \end{eqnarray*}
  $2.$ It suffices to choose for $f \colon X \rightarrow Y$ the identity $1_X \colon X \rightarrow X$.
\end{proof}
\begin{example}
\rm{
In any Goursat variety $\mathcal C$, there is a natural notion of abelian algebra: an algebra $A$ in $\mathcal C$ is abelian if and only if there exists a (necessarily unique) homomorphism $p \colon A \times A \times A \rightarrow A$ satisfying the Mal'tsev identities $p(x,y,y)=x$ and $p(x,x,y)=y$ \cite{Gumm}.
The subcategory ${\mathcal C}_{Ab}$ of abelian algebras is a subvariety of $\mathcal C$: the reflection of an algebra $A$ into ${\mathcal C}_{Ab}$ is simply given by the quotient $\frac{A}{[\nabla_A, \nabla_A]}$ of $A$ by the largest commutator $[\nabla_A, \nabla_A]$ (here $\nabla_A = A \times A$ is the largest congruence on $A$). Accordingly, the closure of a congruence $S$ on $A$ is given by $$\overline{S}= [\nabla_A,\nabla_A] \circ S \circ [\nabla_A,\nabla_A].$$}
\end{example}
It is known that any exact Goursat category $\mathcal C$ is such that the lattice of equivalence relations (=congruences, if $\mathcal C$ in a variety) on any object is modular.
In the following Proposition, we characterize the exact Goursat
categories having the property that the lattice of equivalence relations is distributive,  in terms of a
property of the closure operator.
This is based on a recent
observation due to Bourn \cite{Bourn}, asserting that a Goursat
category is congruence distributive if and only if
$$f(R \wedge S)=f(R) \wedge f(S)$$
for any regular epimorphism $f \colon X \rightarrow Y$ and
equivalence relation $R$ and $S$ on $X$.
\begin{proposition}\label{intersection}
For an exact Goursat
category $ \mathcal C $ the following conditions are equivalent:
\begin{enumerate}
\item the lattice of equivalence relations on any object $X$ in $ \mathcal C $ is
distributive;
\item the closure operator corresponding to any Birkhoff subcategory $ \mathcal L $
of $ \mathcal C $ satisfies the axiom
$$(7) \quad f(\overline{R \wedge S})=\overline{f(R)} \wedge \overline{f(S)}$$
for any regular epimorphism $f \colon X \rightarrow Y$.
\end{enumerate}
\end{proposition}
\begin{proof}
If we assume that $ \mathcal C $ is distributive and $\mathcal L$
a Birkhoff subcategory of $\mathcal C$, we first remark that, for
any equivalence relations $R$ and $S$ on X, one has
\begin{eqnarray*} \overline{R} \wedge  \overline{S}  & = &
(R  \circ \overline{{\Delta}_X}  \circ R) \wedge ( S  \circ
\overline{{\Delta}_X}  \circ
S)  \\
& = & (R  \vee \overline{{\Delta}_X}) \wedge ( S  \vee
\overline{{\Delta}_X}
) \\
 & = &  ( R \wedge S) \vee \overline{\Delta}_X \\
 & = & \overline{R \wedge S .}
  \end{eqnarray*}
  From this, from property $(6)$ of the closure operator and from
  the characterization of the distributive
  Goursat categories recalled above it
  follows that
$$ f(\overline{R \wedge S})= \overline{ f(R \wedge S)} =  \overline{ f(R) \wedge f(S)}
= \overline{ f(R)} \wedge \overline{ f(S)}.$$ Conversely, it
suffices to choose as Birkhoff subcategory $ \mathcal L $ of  $
\mathcal C $ the category  $ \mathcal C $ itself. It follows that,
obviously, $\overline{R} = R$, in that case. Accordingly, by
applying axiom $(7)$, one has that
$$f(R) \wedge  f(S) =\overline{f(R)} \wedge \overline{f(S)}=
f(\overline{R \wedge S})=  f(R \wedge S).$$
\end{proof}
\begin{examples}
\rm{
\begin{enumerate}
\item Recall that a Heyting algebra can be defined as a distributive lattice $(L, \wedge, \vee)$ with a top element $1$ and a bottom element $0$ equipped with
an additional binary operation $\Rightarrow$ satisfying the identities
 $$(a \vee ( (b\Rightarrow a) \wedge b))= a \quad {\rm and} \quad a= (a \wedge (b \Rightarrow (a \wedge b))).$$
If we denote by $\mathbf{Heyting}$ the variety of Heyting algebras, it is well known that it is a Mal'tsev congruence distributive variety (see for instance Example $2.9.16$ in \cite{BB}).
Any subvariety of $\mathbf{Heyting}$ will then fall under the scope of Proposition \ref{intersection}: it is the case, in particular, for the variety $\mathbf{Boole}$ of boolean algebras, which is determined by the additional identity $$((a \Rightarrow 0) \Rightarrow 0) = a.$$
\item A commutative unitary ring $R$ is von Neumann regular if for any element $a$ in $R$ there exists an element $a^*$ such that $$(\alpha )\qquad  a \cdot {a^*}^2= {a^*} \quad {\rm and } \quad a^2 \cdot {a^*}=a.$$ Such an element ${a^*}$ is necessarily unique: this allows one to conclude that the category $\mathbf{CVNReg}$ of commutative von Neumann regular rings is a variety of universal algebras, whose theory is obtained by adding to the theory of unitary commutative rings a unary operation ${( \,)}^*$ satisfying the axioms $(\alpha)$ above. The variety $\mathbf{CVNReg}$ is clearly $2$-permutable, since a group operation is present in its theory; it is also congruence distributive, as shown in Example $2.9.15$ in \cite{BB}. A natural subvariety of $\mathbf{CVNReg}$ to which Proposition \ref{intersection} applies is represented here by the variety $\mathbf{BRng}$ of unitary boolean rings, which is determined by the identity $x^2=x$ (as in the previous example, this subvariety is again the variety $\mathbf{Boole}$ of Boolean algebras, see, for example, \cite{Be}).
\end{enumerate}}
\end{examples}


\begin{thebibliography}{99}
\bibitem{AW} A. Arhangels'skii and R.Wiegandt, \emph{Connectedness and disconnectedness in topology}, Gen. Topology Appl. 5, 1975, 9-33.
\bibitem{Ba} M. Barr, \emph{Exact  Categories},
                 Lect. Notes Math. {236}, Springer-Verlag, 1971, 1-120.
\bibitem{Be} J. Bell and M. Machover, \emph{A course in mathematical Logic},
    North Holland, 1977.
\bibitem{BB} F. Borceux and D. Bourn, \emph{Mal'cev, Protomodular, Homological and Semi-Abelian Categories}, Math. and its Appl. Vol. 566, 2004.
    \bibitem{BC1} F. Borceux and M.M. Clementino, \emph{Topological semi-abelian algebras}, Advances Math., 190, 2005, 425-453.
\bibitem{BC} F. Borceux and M.M. Clementino, \emph{Topological protomodular algebras}, Topology and its Appl., 153, 2006, 3085-3100.
\bibitem{Bourn} D. Bourn, \emph{Congruence distributivity in Goursat and Mal'cev
categories}, Appl. Categ. Structures, 13, 2005, 101-111.
\bibitem{BG} D. Bourn and M. Gran, \emph{Torsion theories in homological categories}, J. Algebra, 305, 2006, 18-47.
\bibitem{CKP}  A. Carboni, G.M. Kelly and M.C. Pedicchio,
     \emph{Some remarks on Maltsev and Goursat categories,}
     Appl. Categ. Struct., {1}, 1993, 385-421.
 \bibitem{CLP} A. Carboni, J. Lambek and M.C. Pedicchio, \emph{Diagram chasing in Mal'cev categories}, J. Pure Appl. Algebra, 69, 1991, 271-284.
 \bibitem{DT} D. Dikranjan and W. Tholen, {\em Categorical structure of closure operator. With applications to topology, algebra and discrete mathematics.} Math. and its Appl. 346, Kluwer Academic Publishers, 1995.
\bibitem{Gumm}  H. P. Gumm, \emph{ Geometrical Methods in Congruence Modular Varieties}, Mem. Amer. Math. Soc. {45}, 286, 1983.
\bibitem{HM} J. Hagemann and A. Mitschke, \emph{On $n$-permutable varieties},
Algebra Univ. 3, 1973, 8-12.
\bibitem{Hus} M. Hu\u sek, \emph{Productivity properties of topological groups},
    Topology Appl. 44, 1992, 189-196.
\bibitem{JK} {G. Janelidze and G.M. Kelly}, \emph{Galois theory and a general
notion of central extension}, {J. Pure Appl. Algebra}, {97}, 1994, 135-161.
\bibitem{Jo} {P.T. Johnstone}, \emph{A note on the semiabelian variety of Heyting semilattices}, in Galois Theory, Hopf Algebras and Semiabelian Categories, The Fields Institute Communications Series, American Mathematical Society, vol. 43, 2004, 317-318.
\bibitem{JP} P.T. Johnstone and M.C. Pedicchio, \emph{Remarks on continuous Mal'cev algebras},
Rend. Istit. Mat. Univ. Trieste 25, no. 1-2, 1994, 277-297.
\bibitem{S} {J.D.H. Smith}, {\em Mal'cev Varieties,} {Lect. Notes Math. 554}, Springer-Verlag, 1976.
\end{thebibliography}
\end{document}